\documentclass{svjour2}                    
\smartqed  
\usepackage{graphicx}
\usepackage{amsmath}
\usepackage{amsfonts}
\usepackage{color}
\usepackage{bm}

\newcommand{\bx}{{\bf X}}

\newcommand{\bv}{{\bf v}}

\newcommand{\be}{\begin{equation}}

\newcommand{\ee}{\end{equation}}
\def\uk{{\bm u}({\bm k},t)}

\def\ukp{u^+({\bm k},t)}
\def\ukm{u^-({\bm k},t)}
\def\bk{\bm k}
\def\bp{\bm p}

\def\bq{\bm q}
\def\bh{\bm h}
\def\bv{\bm v}
\def\bx{\bm x}

\def\hm{{\bm h}^-({\bm k})}
\def\hp{{\bm h}^+({\bm k})}

\newcommand{\cP}{{\cal P}}
\newcommand{\cp}{{\cal P}}

\newcommand{\cH}{{\cal H}}

\begin{document}

\title{On the global regularity  of a helical-decimated version of the 3D Navier-Stokes equations}


\author{Luca Biferale         \and
        Edriss S. Titi 
}


\institute{L. Biferale \at
              Dept. of  Physics and INFN, University of Rome `Tor Vergata', Rome, Italy \\
                          \email{biferale@roma2.infn.it}           
           \and
E.S. Titi          \at
Department  of Computer Science and Applied Mathematics,
               Weizmann Institute of Science, Rehovot 76100, Israel. Also Departments of Mathematics and of Mechanical and Aerospace Engineering, University of California, Irvine, CA 92697,USA. \email etiti@math.uci.edu \& \email{edriss.titi@weizmann.ac.il}
}

\date{Received: March 2,2013 / Accepted: date}

\maketitle

\begin{abstract}
We study the global regularity, for all time and all initial data in $H^{1/2}$, of a recently introduced decimated version of the incompressible
3D Navier-Stokes (dNS) equations. The
model is based on a projection of the dynamical evolution of Navier-Stokes (NS) equations into
the subspace where helicity (the $L^2-$scalar product of velocity and vorticity) is sign-definite. The presence of a second  (beside energy) sign-definite inviscid conserved quadratic quantity, which is equivalent to the $H^{1/2}-$Sobolev norm,
allows us to demonstrate global existence and uniqueness,  of space-periodic solutions,  together with continuity with respect to the initial conditions, for this decimated 3D model. This is achieved thanks to the  establishment of  two new estimates, for this 3D model, which show that the  $H^{1/2}$ and the  time average of the square of the  $H^{3/2}$ norms of the
velocity field remain finite.
Such two additional  bounds are known, in the spirit of the work of H. Fujita and T. Kato \cite{kato1,kato2}, to be sufficient for showing well-posedness for the 3D NS equations. Furthermore, they are directly linked to the helicity evolution for the dNS model, and therefore with a clear physical meaning and consequences.

\end{abstract}

\section{Introduction}
\label{intro}
The problem of whether the three-dimensional incompressible Navier-Stokes equations possess global regular
 and unique solutions is one of the most challenging unsolved problems in applied analysis  \cite{ref1,ref2,ref1_2d,LADY03}. In the two-dimensional case, global existence and uniqueness
 of weak and strong solutions is a well-established result; and in three dimensions, it is possible to show existence of weak solutions,  globally in times, but we do not know how to establish their uniqueness for general initial conditions (see, e.g., \cite{ref1,ref2,ref1_2d} and references therein).
Due to the notorious difficulties to attack the original  problem in its whole generality, many sufficient conditions and criteria have been provided for establishing global regularity (for the most recent development in this regard see, e.g., \cite{berselli2002,cao2008,Cao-Titi,sverak,sverak2,kukavica2006,SESV} and references therein). Moreover,  it may also be important and instructive to start from a {\it decimated} version of it, constraining the evolution of the Navier-Stokes velocity field in sub-manifolds with some given
 physical and dynamical properties. For example, it is well-known that viscous flows with certain physical symmetries, such as axi-symmetric flows \cite{LadyPaper}
and helical flows \cite{MTL},  have global regularity despite the fact that they have nontrivial 3D vortex stretching terms.   Two new examples have recently appeared in the literature. In the first case \cite{frisch2012}, the evolution of the Navier-Stokes equations is projected on a Fractal-Fourier space. In the second case \cite{biferale2012}, the velocity field is first exactly decomposed for each Fourier mode in two components, carrying positive and negative helicity respectively, and then evolved by keeping only one of the two, such as to make  also the helicity  invariant sign-definite (i.e. positive). Both approaches have been proposed in order to understand key phenomenological properties of  turbulence as a function of its total number of degrees-of-freedom (for the Fractal-Fourier decimation) or as a function of the physical character of the non-linear interactions (for the helical-Fourier decimation; this is to distinguish from the above mentioned flows with physical helical symmetry \cite{Ettinger,MTL}). Clearly, it is tempting to ask also questions concerning regularities of their solutions. By investigating this model we hope to  learn about possible crucial physical features of the ``long sought'' singular solutions of the original Navier-Stokes equations that are eventually absent in our decimated version and in models with physical spatial symmetries. A similar approach has recently been proposed in \cite{grauer2012}
where the existence of global regularity for solutions of a modified Navier-Stokes equation has been proven. Let us notice however, that at difference from the previous cases, the equations studied in \cite{grauer2012}  cannot be seen as a {\it decimation} of the original Navier-Stokes case.

In this paper we will concentrate on the regularity properties of the Navier-Stokes equations under helical-Fourier decimation only. Without loss of generality, we will study the case where the projection is made into the sub-manifold where helicity is positive-definite. It is worth mentioning that in this situation  the positive-definite helicity is equivalent to the square of the $H^{1/2}-$Sobolev norm.

In such a case, we are able to show that the dNS system possesses global unique solutions which depend continuously on the initial data.
Such a result is obtained thanks to two additional estimates available only when the dynamics is projected on a manifold with a sign-definite helicity.  From a phenomenological physical point of view,
this  dNS equations also show a {\it non-typical}  inverse-energy cascade, i.e. when forced on an
intermediate range in Fourier space, energy tends to go to larger and larger scales and helicity is transported
toward small-scales \cite{biferale2012}.
This is the so-called split-cascades regime  typical of 2D Navier-Stokes flows - where enstrophy flux replaces the helicity (which is a positive-definite quantity here) flux.
Such a phenomenology has also been observed in real 3D Navier-Stokes equations under fast rotation and with helical forcing \cite{pouquet_dns,smith96}.

The article is organized as follows. In section (\ref{sec:he}) we
review the main properties of the helical-decimated equations. In sec. (\ref{sec:ex})  we briefly discuss the existence of weak solutions to the dNS system following the usual construction based on the Galerkin approximation procedure. In
secs. (\ref{sec:un}-\ref{sec:more})
we establish  two additional global uniform  bounds
on the velocity  field, for initial data with finite positive helicity,  by exploiting the constraints from the inviscid conservation of helicity. Furthermore, we use these additional estimates to prove  the uniqueness of the more regular solutions of the dNS system, i.e. solutions with initial data of finite positive helicity.  Conclusions and a brief discussion on
further perspectives can be found in  sec. (\ref{sec:co}).

\section{The helical-decimated Navier-Stokes equations}
\label{sec:he}
\begin{figure*}
\begin{center}
\includegraphics[width=0.6\textwidth]{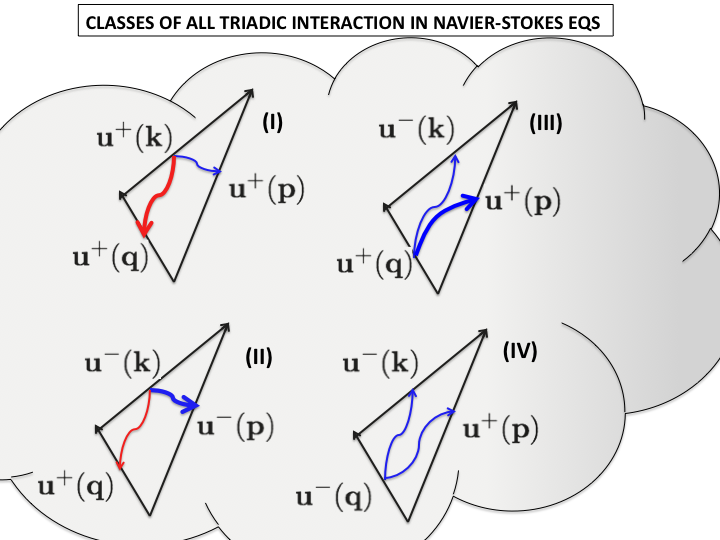}
  \caption{Under helical decomposition the three wavenumbers interaction of the non-linear NS
equations are decomposed in four classes, depending on the relative helicity signs.
In \cite{waleffe} a simple dynamical argument is given supporting the fact that triads of classes
(III) and (IV) mainly transfer energy toward small scales (high wave-numbers), i.e.
they have the usual direct cascade, triads of class (I) enjoys an inverse energy cascade,
while class (II) is mixed. In the figure this is summarized by red arrows denoting a
backward energy transfer and by blue arrows for forward energy transfer.
In \cite{biferale2012} a direct numerical integration at
 high resolution of the dNS
with only  triads of class (I) showed that a stationary turbulent inverse energy cascade is
indeed established.}
\label{fig:waleffe}
\end{center}
\end{figure*}

The starting point of our analysis is the well-known helical-Fourier  decomposition \cite{waleffe} of the velocity field  $\bv(\bx,t)$,
expanded in Fourier series, $\uk$. Here we consider spatially periodic flows with zero spatial mean; consequently the zero wavenumber Fourier mode ${\bm u}({\bm 0},t)=\bm 0$.  Being divergence-free, each velocity component
 in Fourier space has only two degrees of freedom.
The idea is to define the two independent degrees of freedom by a projections on two orthonormal vectors, each one  bringing a definite sign of helicity:
\begin{equation}
  \uk  = \ukp \hp + \ukm \hm
\end{equation}
where $ {\bm h}^\pm(\bk)$ are the eigenvectors of the curl operator $i {\bm
  k} \times {\bm h}^\pm(\bk) = \pm k {\bm h}^\pm(\bk)$.
 In particular,  one can always  choose $ {\bm h}^\pm(\bk) =\hat{\bm \mu}(\bk) \times
\hat{\bm k} \pm i \hat{\bm \mu}$, where $ \hat{\bm \mu}$ is an
arbitrary unit vector orthogonal to ${\bm k}$ which satisfies the relation
$\hat{\bm \mu}({\bm k}) = - \hat{\bm \mu}(-{\bm k})$ (necessary to
ensure the reality of the velocity field).  Such requirement is
satisfied, e.g., by the choice $\hat{\bm \mu}(\bk) = {\bm z} \times {\bm k}
/|| {\bm z} \times {\bm k} ||$, with ${\bm z}$ an arbitrary vector.
 We then have for energy, $E(t) = \int d^3 x \, |\bv(\bx,t)|^2$, and helicity, $\cH(t) = \int d^3 x \, \bv(\bx,t) \cdot \bm{\omega} (\bx,t)$ (here $\bm{\omega}(\bx,t) = \nabla \times \bv(\bx,t)$ is the vorticity):
\begin{equation}
  \begin{cases}
    E(t) = \sum_{\bk} |\ukp|^2 + |\ukm|^2; \label{eq:E}\\
    \cH(t) = \sum_{\bk} k(|\ukp|^2 - |\ukm|^2).
  \end{cases}
\end{equation}
The non-linear term of the NS equations can be exactly
decomposed in 4 independent classes of triadic interactions,
determined by the helical content of the complex amplitudes, ${u}^{s_k}(\bk)$ with $s_k = \pm $
(see \cite{waleffe} and Fig.~\ref{fig:waleffe}).
 Among three generic interacting modes ${u}^{s_k}(\bk),{ u}^{s_p}(\bp), {
  u}^{s_q}(\bq)$, one can identify 8 different helical combinations
$(s_k=\pm,s_p=\pm,s_q=\pm)$.  Among them, only four are independent
because of the symmetry that allows to change all signs of helicity simultaneously.\\
Therefore, the Navier-Stokes equations in the helical-Fourier basis are \cite{waleffe}:
\begin{eqnarray}
(\partial_t + \nu k^2)\overline{u}^{s_k}(\bk,t) =& \\
 -\frac{1}{4} \sum_{\bk+\bp+\bq=0}  \sum_{s_p, s_q} (s_pp-s_q q) &
[\bh^{s_p}(\bp)\times\bh^{s_q}(\bq)\cdot \bh^{s_k}(\bk)]u^{s_p}(\bp,t)u^{s_q}(\bq,t) + \overline{f}^{s_k}(\bk) \nonumber
\end{eqnarray}
where $\nu$ is the viscosity,  ${f}^{s_k}(\bk)$, the external forcing and where 
$\overline{\cdot}$ stands for complex conjugate.
It is important to remark that the non-linear dynamics formally  preserves both
energy and helicity, triad-by-triad, i.e.
$\partial_t (|u^{s_p}|^2+|u^{s_q}|^2+|u^{s_k}|^2) = \partial_t(s_p p|u^{s_p}|^2+s_qq|u^{s_q}|^2+s_kk|u^{s_k}|^2) =0$. As a result, the evolution for  helicity ``formally" becomes:
\begin{eqnarray}
\frac{d}{dt}  \cH(t)  = -\nu \sum_{\bk} k^3 (|\ukp|^2 - |\ukm|^2).+ &\nonumber \\
 2 Re \sum_{\bk} k ( u^+(\bk,t)\overline{f}^+(\bk)-u^-(\bk,t)\overline{f}^-(\bk)).
\label{eq:hfull}&
\end{eqnarray}
Among the 4 possible different choices of  helicity signs interacting in the nonlinear
term we will concentrate here on what happens
when only one class is present. In particular, we will ask what happens when only interactions having the same sign of helicity for all the three wavenumber in each triad, i.e. helicity becomes sign-definite. In order to do that  we define the projector on
positive/negative helicity states as
\begin{equation}
  \label{eq:proj}
  {\cp}^\pm (\bk)\equiv
  \frac {{\bm h}^\pm (\bk) \otimes \overline{{\bm h}^{\pm}}(\bk)}
  {\overline{{\bm h}^\pm}(\bk) \cdot {\bm h}^\pm(\bk)}.
\end{equation}
The projector is self-adjoint and commutes with derivatives.
 The action of this projector on a field in real space is defined via its Fourier
decomposition, for example for projection on positive helicity we have:
\begin{equation}
  \label{eq:projv}
  {\bv}^+(\bx,t) \equiv  {\mathcal P}^+ \bv(\bx,t) \equiv \sum_{\bk} e^{i{\bm k}\cdot \bx}{\cp}^+(\bk) {\uk} \equiv \sum_{\bk} e^{i{\bm k}\cdot\bx} \ukp \hp .
\end{equation}
We can then consider the dynamical evolution of the dNS:
\begin{equation}
  \label{eq:ns+++}
\begin{cases}
  \partial_t {\bm v^+} =   {\mathcal P}^+(- {\bm v^+} \cdot {\bm \nabla} {\bm v^+} -{\bm \nabla} p^+ )
  +\nu \Delta {\bm v^+} + {\bm f}^+\\
\nabla \cdot {\bm v^+} =0
\end{cases}
\end{equation}
where $p$ is the pressure and ${\bm f}^+$ is the
external forcing, here taken for simplicity time independent and with zero mean such as we can restrict the analysis of the solutions of
(\ref{eq:ns+++}) in the space where $\int {\bm v^+} dx =0$. Following the tradition  we denote by $\dot{X}$ the linear subspace of all functions in $X$ that have spatial mean zero. Notice, that it is easy to realize that if the initial velocity configuration is chosen with only positive helicity components, i.e.
${\bv}^-(\cdot,t=0)=0$, the dynamical evolution of (\ref{eq:ns+++}) preserves this property
 for all times. \\

\section{Global existence of weak solutions with initial data in $\dot{L}^2$}
\label{sec:ex}
The notion of weak solution for the dNS system (\ref{eq:ns+++}) is similar to that of the 3D NS equations \cite{ref1,ref1_2d}. Specifically, for a given    initial data $\bm v^+(\cdot,t=0)=\bm v^+_0(\cdot)\in \dot{L}^2$, a weak solution of (\ref{eq:ns+++}) in the interval $[0,T]$ is a  divergence-free vector field $\bm v^+ (\bm x, t)$, which  belongs to   $ C_{\hbox{weak}}([0,T], \dot{L}^2) \bigcap L^2([0,T], \dot{H}^1)$, and for which $\partial_t{\bm v^+}\in L^{5/4}([0,T], \dot{H}^{-1})$; such that  the first equation of (\ref{eq:ns+++}) holds in $L^{5/4}([0,T], \dot{H}^{-1})$.
It is not difficult to prove the existence of weak solution for the dNS. It is enough to proceed as for the standard NS case, introducing
a Galerkin projection on all wavenumbers smaller than a given value, $N$:
\be
\bv_N(\bx,t) = \cP_N \bv(\bx,t) = \sum_{|\bk| <N} e^{i{\bm k}\cdot \bx} {\uk},
\ee
  and then controlling the limit of the Galerkin approximates for $ N \rightarrow \infty$.
The only additional detail consists in checking that the Galerkin projector, $\cP_N$ commutes
with the helical-Fourier projector $\cP^+$, which is of course the case. We can then define a Galerkin-helical-Fourier
projected field, evolving according to the finite system of ordinary differential equations:
\begin{equation}
  \label{eq:ghf}
\begin{cases}
  \partial_t {\bm v^+_N} =   {\mathcal P}^+ \cP_N (- {\bm v^+_N} \cdot {\bm \nabla} {\bm v^+_N} -{\bm \nabla} p^+_N )
  +\nu \Delta {\bm v^+_N} + {\bm f}_N^+\\
\nabla \cdot {\bm v^+_N} =0.
\end{cases}
\end{equation}
Since the nonlinearity in the above ODE system is quadratic it has a short time existence and uniqueness solution for all initial data. As in the NS equations case one can show that the $L^2-$norm of the solutions to the above system remain bounded and hence it has global existence and uniqueness. Following the same steps   as for the standard original
NS case  \cite{ref1,ref2} starting from any finite $N$
and then passing to the limit, modulo a subsequence  $N_j \rightarrow \infty$, using the relevant compactness theorem of Aubin \cite{ref1}.

\section{More regular solutions with initial data in $\dot{H}^{1/2}$}
\label{sec:more}
In this section, we are going to use the fact that for the inviscid and unforced  dNS (\ref{eq:ghf}) system  the helicity is formally conserved, and that it is positive-definite quadratic quantity, which is equivalent to the square of the $\dot{H}^{1/2}-$ Sobolev norm. Therefore, obtaining uniform (in time) bounds on the helicity enables us to prove the existence of solutions with a higher degree of regularity, provided the initial data is in $\dot{H}^{1/2}$. Furthermore, this additional regularity will allow us to prove, in the next section,  the uniqueness of these regular solutions within the class of weak solutions.  Let us observe that all  the  estimates that follow are formal, but can be rigorously justified
 by obtaining them first for the corresponding solutions of the     Galerkin approximating system (\ref{eq:ghf}), and then passing to the limit, modulo subsequences,  with $N \to \infty$. Furthermore, it is worth mentioning that similar ideas and estimates can be found in \cite{kato1,kato2} in the study of short time existence and uniqueness of the three-dimensional NS equations with initial data in $H^{1/2}$. The advantage of  system (\ref{eq:ns+++}) over the NS equations is that the $H^{1/2}$ remains finite, which allows to extend the short time existence argument to prove global regularity for all time and all initial data in $H^{1/2}$. Indeed, the fact that  helicity is  a (positive-definite) inviscid
invariant of (\ref{eq:ns+++}) can be readily
verified, formally, by its evolution written in terms of the Fourier components (\ref{eq:hfull}), which after decimation reads:
\be
\label{eq:helicity}
\frac{1}{2}\frac{d}{dt} \sum_{\bk} k |u^+(\bk,t)|^2  + \nu \sum_{\bk} k^3 |u^+(\bk,t)|^2 = Re \sum_{\bk} k u^+(\bk,t)\overline{f}^+(\bk),
\ee
where $k = |\bk|$, and $u^+(\bk,t) = {\cp}^+(\bk) {\uk} $. Furthermore,
 the right-hand side (RHS) of  (\ref{eq:helicity}) can be bounded by:
\be
\label{eq:helicity2}
 RHS \le \sum_{\bk}|k^{\frac{3}{2}} u^+(\bk,t)|\, |f^+(\bk) \,k^{-\frac{1}{2}}| \le (\sum_{\bk} k^3 |u^+(\bk,t)|^2)^{\frac{1}{2}}(\sum_{\bk} |f^+(\bk)|^2k^{-1})^{\frac{1}{2}},
\ee
where the last step is obtained by using Cauchy-Schwarz inequality. By multiplying and dividing
by $\sqrt{\nu}$ and   using Young's inequality  we can further bound it by:
\begin{eqnarray}
(\nu \sum_{\bk} k^3 |u^+(\bk,t)|^2)^{\frac{1}{2}}(\sum_{\bk} |f^+(\bk)|^2k^{-1}\nu^{-1})^{\frac{1}{2}} \le\\
 \frac{\nu}{2} \sum_{\bk} k^3 |u^+(\bk,t)|^2
+\frac{1}{2 \nu} \sum_{\bk} |f^+(\bk)|^2k^{-1}.
\end{eqnarray}
Plugging back the above estimate in (\ref{eq:helicity}) we get:
\be
\label{eq:helicity3}
\frac{1}{2}\frac{d}{dt} \sum_{\bk} k |u^+(\bk,t)|^2  + \frac{\nu}{2} \sum_{\bk} k^3 |u^+(\bk,t)|^2 \le \frac{1}{2 \nu} \sum_{\bk} |f^+(\bk)|^2k^{-1}.
\ee
Observe that for  every  $\phi \in \dot H^s$ (i.e.  $\phi$ belongs to the Sobolev space $H^s$  of   periodic functions  with zero average, $\int \phi(x) \, dx =0$) the
 $H^s-$Sobolev norm, for $s\ge 0$, is equivalently defined by
\be
||\phi||_{H^s} = (\sum_{\bk} k^{2s} |\phi(k)|^2)^{1/2}.
\ee
Thanks to the above observation we can rewrite (\ref{eq:helicity3})
in terms of the $H^{\frac{1}{2}}$ and $H^{\frac{3}{2}}$ norms as:
\be
\label{eq:helicity4}
\frac{1}{2}\frac{d}{dt} ||v^+||_{H^{\frac{1}{2}}}^2  + \frac{\nu}{2} ||v^+||_{H^{\frac{3}{2}}}^2 \le \frac{1}{2 \nu} \sum_{\bk} |f^+(\bk)|^2k^{-1}.
\ee
Expression (\ref{eq:helicity4}) can be rewritten in
 its equivalent integral form as:
\be
\label{eq:helicity5}
||v^+(t)||_{H^{\frac{1}{2}}}^2  + \int_0^t dt' \frac{\nu}{2} ||v^+(t')||_{H^{\frac{3}{2}}}^2 \le ||v^+(0)||_{H^{\frac{1}{2}}}^2+\int_0^t dt' \frac{1}{2 \nu} \sum_{\bk} |f^+(\bk,t')|^2k^{-1}.
\ee
The above inequality  implies   that  there is a weak solution of the dNS (\ref{eq:ns+++}) which satisfies:
\begin{equation}\label{Reg}
v^+  \in   L^{\infty}([0,T]; \dot{H}^{\frac{1}{2}}); \qquad
 v^+  \in  L^{2}([0,T]; \dot{H}^{\frac{3}{2}}),
\end{equation}
provided that ${\bm f} \in L^{2}([0,T];  H^{-\frac{1}{2}})$ and ${\bm v^+}(0) \in \dot H^{\frac{1}{2}}$. In the next section we will show that this more regular solution, that satisfies (\ref{Reg}), is unique within the class of weak solutions.

\section{Uniqueness of more regular solution}
\label{sec:un}
In order to prove the uniqueness of the more regular solutions with initial data in $\dot{H}^{1/2}$, within the class of weak solutions, let us first
 rewrite the dNS equations (\ref{eq:ns+++}) in a symbolic form:
\begin{equation}
  \label{eq:ns}
  \partial_t {\bm v^+} -\nu \Delta {\bm v^+} =
- {\bm B}_{\cH}({\bm v^+},{\bm v^+})  + {\bm f}^+\\
\end{equation}
where $ {\bm B}_{\cH}({\bm a},{\bm b}) = \cP^+ [({\bm a} {\bm \nabla})\,{\bm  b} - \Delta^{-1}{\bm \nabla} {\bm \nabla}  {\bm a} {\bm b})]$\\
Below, we give a formal proof of uniqueness which can be made fully rigorous by working  with the
integrated -in time- version of (\ref{eq:ns2}) and then following exactly the same argument as in Serrin \cite{serrin} (see also \cite{bardos}). It is worth mentioning that similar ideas and estimates can be found in  \cite{kato1,kato2} in the study of short time well-posedness of three-dimensional NS equations with ${H}^{1/2}$ initial data.
Let us begin by considering two different weak solutions,
${\bm v^+_1}$ and ${\bm v^+_2}$, with the same initial data in $\dot{H}^{1/2}$,  and where ${\bm v^+_1}$ satisfies  (\ref{Reg}). Their difference, ${\bm w^+} = {\bm v}^+_1- {\bm v}^+_2$, satisfies the evolution equation:
\begin{equation}
  \label{eq:ns2}
  \partial_t {\bm w^+} - \nu \Delta {\bm w^+} =  - {\bm B}_{\cH}({\bm v^+_1},{\bm w^+})  - {\bm B}_{\cH}({\bm w^+},{\bm v^+_1})  + {\bm B}_{\cH}({\bm w^+},{\bm w^+}). \\
\end{equation}
Next, we study the time evolution of $|| w^+ ||_{L^2}^2 = \int d^3x | {\bm w^+}|^2 $, by formally multiplying (\ref{eq:ns2}) and integrating over the spatial domain.
Note that  $\cP^+$ is a self-adjoint operator, and that
$\int d^3 x \, (\cP^+ {\bm a^+})\cdot {\bm b^+} = \int d^3 x \, {\bm a^+}\cdot  \cP^+ {\bm b^+} = \int d^3 x\, {\bm a^+}\cdot {\bm b^+}$.
By incompressibility, and after integration by parts and using the periodic boundary conditions,  only the second term in the RHS of (\ref{eq:ns2}) survives and we formally obtain:
\begin{equation}
  \label{eq:ns1}
\frac{1}{2}  \frac{d}{dt} || w^+||_{L^2}^2 +\nu  ||\nabla w^+||_{L^2}^2 = - \int d^3x  [({\bm w^+}\cdot {\bm  \nabla}) {\bm  v^+_1}]\cdot {\bm w^+}.
\end{equation}
The RHS of (\ref{eq:ns1}) can be further integrated by part, using once again the incompressibility and the periodic boundary conditions,  to get
$\int d^3x  [({\bm w^+}\cdot {\bm  \nabla}) {\bm  v^+_1}]\cdot {\bm w^+} = - \int d^3x  [({\bm w^+}\cdot {\bm  \nabla}) {\bm  w^+}]\cdot {\bm v^+_1}$; and by applying H\"{o}lder's inequality one obtains the following bound:
\be
\label{eq:estimate}
|\int d^3x  [({\bm w^+}\cdot {\bm  \nabla})\, {\bm  w^+}]\cdot {\bm v^+_1}|  \le || w^+||_{L^3} || \nabla w^+ ||_{L^2} ||v^+_1||_{L^6}
\ee
 Applying the three-dimensional  Sobolev embedding inequalities:
\[   ||\cdot||_{L^3}  \le const. ||\cdot||_{H^{1/2}} \quad
||\cdot||_{L^6}  \le const. ||\cdot||_{H^{1}}
\]
to (\ref{eq:estimate})  we obtain:
\be
|\int d^3x  [({\bm w^+}\cdot {\bm  \nabla})\, {\bm  w^+}]\cdot {\bm v^+_1}|  \le  const. || w^+||_{H^{\frac{1}{2}}} || \nabla w^+ ||_{L^2} ||v^+_1||_{H^1}.
\ee
We now use the interpolation inequality $||\cdot ||_{H^s} \le  || \cdot||^{\theta}_{H^{s_1}} || \cdot||^{(1-\theta)}_{H^{s_2}}$, with
$s = \theta s_1 + (1-\theta) s_2$ with $\theta = 1/2, s_1=0, s_2=1$ to get:
\be
||w^+||_{H^{\frac{1}{2}}} \le  ||w^+||^{\frac{1}{2}}_{H^0} ||w^+||^{\frac{1}{2}}_{H^1} =  ||w^+||^{\frac{1}{2}}_{L^2} ||\nabla w^+||^{\frac{1}{2}}_{L^2};
\ee
and with   $\theta = 1/2, s_1=1/2, s_2=3/2$ to obtain:
\be
||v^+_1||_{H^1} \le ||v^+_1||^{\frac{1}{2}}_{H^{\frac{1}{2}}} ||v^+_1||^{\frac{1}{2}}_{H^{\frac{3}{2}}}.
\ee
Putting all together we can bound the LHS of  (\ref{eq:estimate}) as:
\be
\label{eq:estimate2}
\le const. ||\nabla w^+||^{\frac{3}{2}}_{L^2} ||w^+||^{\frac{1}{2}}_{L^2}   ||v^+_1||^{\frac{1}{2}}_{H^{\frac{1}{2}}} ||v^+_1||^{\frac{1}{2}}_{H^{\frac{3}{2}}}.
\ee
Let us now multiply and divide by  $\nu^{\frac{3}{4}}$:
\be
\label{eq:estimate3}
\le  const. [\nu^{\frac{3}{4}} ||\nabla w^+||^{\frac{3}{2}}_{L^2}] \frac{||w^+||^{\frac{1}{2}}_{L^2}   ||v^+_1||^{\frac{1}{2}}_{H^{\frac{1}{2}}} ||v^+_1||^{\frac{1}{2}}_{H^{\frac{3}{2}}}}{\nu^{3/4}};
\ee
and apply Young's inequality with the weights $3/4$ and $1/4$:
\be
\label{eq:estimate4}
\le \frac{3}{4} \nu ||\nabla w^+||^2_{L^2}] + const. \frac{||w^+||^{2}_{L^2}   ||v^+_1||^{2}_{H^{\frac{1}{2}}} ||v^+_1||^{2}_{H^{\frac{3}{2}}}}{\nu^3}.
\ee
Going back to equation (\ref{eq:ns1}) we can then bring the first term of the above relation on the LHS and get:
\begin{equation}
  \label{eq:ns3}
\frac{1}{2}  \frac{d}{dt} || w^+||_{L^2}^2 +\frac{1}{4} \nu  ||\nabla w^+||_{L^2}^2 \le K \frac{||w^+||^{2}_{L^2} ||v^+_1||^{2}_{H^{\frac{3}{2}}}}{\nu^3},
\end{equation}
where we have used the assumption that  $v_1^+  \in   L^{\infty}([0,T]; \dot{H}^{\frac{1}{2}})$, for every $T>0$, since it satisfies (\ref{Reg}). By virtue of
Gronwall's inequality we finally arrive to:
\be
||w^+(t) ||_{L^2}^2 \le ||w^+(0) ||_{L^2}^2  \exp{\Big (K \nu^{-3} \int_0^t dt' ||v^+_1(t')||^{2}_{H^{\frac{3}{2}}}\Big )}
\ee
which proves uniqueness of the solution thanks to the fact that we have
$  v_1^+ \in L^2([0,T]; \dot H^{\frac{3}{2}})$ (since it satisfies (\ref{Reg}). The above also implies  the Lipschitz continuity with respect to the initial data due to the $||w^+(0) ||_{L^2}^2 $
prefactor.

\section{Conclusions and further perspectives}
\label{sec:co}
Inviscid invariants of the NS equations are crucial in determining the
physics of turbulence \cite{frisch}.
In 2D turbulence the presence of two positive-definite quadratic quantities,
energy and enstrophy, leads to a split cascade: energy flows toward large scales and enstrophy to small scales
\cite{boffi_2d,boffi_2da,2dreview}.  The fluid equations possess two inviscid invariants also
in 3D: energy and helicity. At variance with energy, helicity is in general not sign-definite, and therefore it is not a `coercitive' quantity.
In principle,  this allows for a simultaneous forward  transfer of energy and helicity, as confirmed by the results
of two-point closures~\cite{brissaud,waleffe} and direct numerical
simulations~\cite{chen,chen2,pouquet_dns}.  Nevertheless, a reversal of the
flux of energy has been observed in geophysical flows subject to Earth's
rotation \cite{pouquet_dns,smith96} if the forcing is chiral, injecting
a non-vanishing helicity in the system. For fast rotating turbulence in the inverse energy cascade regime, helicity is observed to still flow toward small scales: we have
a split-cascade scenario as for the energy-enstrophy case of 2D turbulence.
Further, the link between  the existence of intermittent burst in the energy cascade and a
 ``local'' helicity blocking mechanism has been
proposed~\cite{biferale_hel}.

In this paper we have shown that helicity plays also a very peculiar role in the determining the `regularity' properties of the velocity field. If
the dynamics is restricted to the sub-set of  modes with a well definite sign of helicity (i.e. positive) then the flow admits unique global weak solutions that depend continuously on the initial data. In our system the vortex stretching mechanism does not vanish, even
 if it is certainly depleted because of the local 'Beltrami' condition imposed at each wave number by our decimation rule.
 Such a result leads to speculation that possible `singular' solution of the original NS equations should then be searched elsewhere, i.e. in the three classes of non-linear interactions connecting triads with different helicity components which are killed in our decimated system.

\begin{acknowledgements}
L.B. would like to thank U. Frisch for stimulating the interest on the problem and V. Zheligovsky for a very careful reading of the manuscript. L.B. also acknowledges the useful discussions with
C. Doering, U. Frisch, A. Porretta and V. Sverak. The work of E.S.T. is  supported in part by the NSF grants no. DMS-1009950, DMS-1109640,  and DMS-1109645. E.S.T. also acknowledges the support of the  Minerva Stiftung/Foundation.
\end{acknowledgements}


\begin{thebibliography}{}
%
%
\bibitem{bardos} C. Bardos, M.C.L. Filho, D. Niu, H.J.N. Lopes and E.S. Titi ``Stability of two-dimensional viscous incompressible flows under three-dimensioanal perturbations and inviscid symmetry breaking'' arXiv:1201.2742v2 (2012).

\bibitem{berselli2002} L. C. Berselli and G. Galdi.
``Regularity criteria involving the pressure for the weak solutions to the Navier-Stokes equations'' Proc. Amer. Math. Soc. {\bf 130} 3583 (2002).



\bibitem{biferale_hel}L. Biferale. ``Shell models of energy cascade in turbulence''  {\it Ann. Rev. Fluid Mech.} {\bf 35} 441 (2003).

\bibitem{biferale2012} L. Biferale, S. Musacchio and F. Toschi. Phys. Rev. ``Inverse energy cascade in three-dimensional isotropic turbulence'' Lett. {\bf 108} 164501 (2012).



\bibitem{boffi_2d} G. Boffetta \& R.E. Ecke. ``Two-dimensional turbulence''          {\it Annu. Rev. Fluid Mech.} 44, 427 (2012);

\bibitem{boffi_2da}     G. Boffetta and S. Musacchio. ``Evidence for the double cascade scenario in two-dimensional turbulence''
        {\it Phys. Rev. E} {\bf 82} 016307 (2010).


\bibitem{brissaud} A. Brissaud, U. Frisch, J. Leorat, M. Lesieur and
  M. Mazure. ``Helicity cascades in fully developed isotropic turbulence'' {\it Phys. Fluids} {\bf 16} 1366 (1973).


\bibitem{cao2008} C.~Cao and E.S.~Titi. ``Regularity criteria for the three--Dimensional Navier-Stokes equations''
{\it Indiana Univ. Math. J.}, {\bf 57}  2643 (2008).

\bibitem{Cao-Titi} C.~Cao and E.S.~Titi. ``Global  regularity criterion  for the $3D$  Navier--Stokes equations involving one entry of the
velocity gradient tensor''  {\it Archive of Rational
 Mechanics \& Analysis},  {\bf 202}  (2011),  919--932.


\bibitem{chen} Q. Chen, S. Chen \& G.L. Eyink. ``The joint cascade of energy and helicity in three-dimensional turbulence'' {\it
    Phys. Fluids} {\bf 15} 361 (2003)
\bibitem{chen2}  Q. Chen, S. Chen, G.L. Eyink.
and D.D. Holm. ``Intermittency in the joint cascade of energy and helicity'' {\it Phys. Rev.Lett} {\bf 90} 214503 (2003).

\bibitem{2dreview} H.J.H. Clercx \& G.J.F. van Heij ``Two-dimensional Navier-Stokes turbulence in bounded domains''
  {\em     Appl. Mech. Rev.} 62, 020802 1-25 (2009)

\bibitem{ref1} P. Contantin and C. Foias. ``Navier-Stokes Equations'' (The University of Chicago Press, 1988)


\bibitem{ref2}
C.R. Doering and J.D. Gibbon. ``Applied Analysis of the Navier-Stokes Equations'' (Cambridge University Press 2001).

\bibitem{sverak} L. Escauriaza, G. Seregin and V. Sverak. ``Backward uniqueness for parabolic equations'' {\em Arch. Ration. Mech. Anal.} {\bf 169} 147 (2003).


\bibitem{Ettinger} B. Ettinger and E.S. Titi. ``Global existence
and uniqueness of weak solutions of 3-D Euler equations with helical
symmetry in the absence of vorticity stretching'' {\em  SIAM, Jour.
Math. Anal. }  {\bf 41(1)} 269 (2009).


\bibitem{frisch} U. Frisch {\it Turbulence: the legacy of
    A.N. Kolmogorov} (Cambridge Univ. Press, 1995).

\bibitem{frisch2012} U. Frisch, A. Pomyalov, I. Procaccia and S.S. Ray.
``Turbulence in Noninteger Dimensions by Fractal Fourier Decimation'' {\em Phys. Rev. Lett.} {\bf 108} 074501 (2012).

\bibitem{kato1} H. Fujita and T. Kato. ``On the Navier-Stokes initial value problem'' {\em Archive for Rational Mech. and Analysis} {\bf 16} 269-315 (1964).


\bibitem{grauer2012} T. Grafke, R. Grauer and T. Sideris. ``Turbulence properties and global regularity of a modified Navier-Stokes equation''
{\em  Physica D} (submitted) (2012).

\bibitem{sverak2} H. Jia and V. Sverak. ``Minimal $L^3-$initial data for potential Navier-Stokes singularities'' arXiv:1201.1592v1 (2012).

\bibitem{kato2} T. Kato and H.  Fujita. ``On the non-stationary Navier-Stokes system'' {\em  Rend. Sem. Mat. Univ. Padova}
 {\bf 32}, 243-260 (1962).

\bibitem{kukavica2006} I. Kukavica and M. Ziane. ``One component regularity for the Navier Stokes equations'' {\em Nonlinearity} {\bf 19} 453 (2006).


\bibitem{ref1_2d}
O. Ladyzhenskaya. ``The Mathematical Theory of Viscous Incompressible Flows'' (New York, Gordon and Breach 1969).


\bibitem{LadyPaper} O.A.~Ladyzhenskaya. ``On unique solvability of three-dimensional
Cauchy problem for the Navier-Stokes equations under the axial symmetry'' {\em Zap. Nauchn. Sem. LOMI} {\bf 7}  155 (1968) (In Russian)


\bibitem{LADY03} O.A.~Ladyzhenskaya. ``The sixth millennium problem: Navier--Stokes equations, existence and smoothness'' (Russian). {\em Uspekhi Mat. Nauk.} {\bf 58(2)}    45 (2003); translation in   {\em  Russian Math.
Surveys}  {\bf 58(2)}  251 (2003).

\bibitem{MTL} A.~Mahalov, E. S.~Titi, and S.~Leibovich. ``Invariant helical subspaces for the Navier-Stokes
equations''  {\em Arch. Rat. Mech. Anal.}   {\bf 112}  193 (1990).


\bibitem{pouquet_dns} P.D. Mininni, A. Pouquet. ``Rotating helical turbulence. I. Global evolution and spectral behavior''  {\it
    Phys. Fluids} {\bf 22} 035105 (2010).


    


\bibitem{SESV} G.~Seregin and V.~Sver\'ak. ``Navier-Stokes equations with lower bounds on the pressure''  {\em Arch.~Rational Mech.~Anal.}
 {\bf 163}  65 (2002).
 
\bibitem{serrin} J. Serrin. ``The initial value problem for the Navier-Stokes equations'' Nonlinear Problems (University of Wisconsin Press, Madison, R. E. Langer edition 1963). 

\bibitem{smith96} L.M. Smith, J.R. Chasnov, F. Waleffe. ``Crossover from two-to three-dimensional turbulence''  {\it Phys. Rev. Lett.}
  {\bf 77}, 2467 (1996).

\bibitem{waleffe} F. Waleffe. ``The nature of triad interactions in homogeneous turbulence'' {\it Phys Fluids A} {\bf 4} 350 (1992).



\end{thebibliography}


\end{document}